\documentclass[12pt,reqno]{article}
\usepackage{amssymb,amscd,amsmath,amsthm,color}
\newtheorem{theorem}{Theorem}[section]

\newtheorem{cor}[theorem]{Corollary}
\newtheorem{prop}[theorem]{Proposition}

\usepackage{enumerate,amssymb}
\theoremstyle{definition}

\newtheorem{example}[theorem]{Example}

\theoremstyle{remark}
\newtheorem{remark}[theorem]{Remark}

\numberwithin{equation}{section}

\def\bC{\mathbb{C}}

\def\bM{\mathbb{M}}

\begin{document}
\baselineskip=15pt

\title{On the Russo-Dye Theorem for positive linear maps}

\author{ Jean-Christophe Bourin{\footnote{This research was supported by the French Investissements  d'Avenir program, project ISITE-BFC (contract ANR-15-IDEX-03).
}} \ and  Eun-Young Lee{\footnote{This research was supported by
Basic Science Research Program through the National Research
Foundation of Korea (NRF) funded by the Ministry of
Education (NRF-2018R1D1A3B07043682)}  }}

\date{ }

\maketitle

\vskip 10pt\noindent
{\small
{\bf Abstract.}  {We revisit a classical result, the Russo-Dye Theorem, stating that every  positive linear map attains its norm at the identity.
}
\vskip 5pt\noindent
{\it Keywords.} Matrix inequalities,  Positive linear maps,  Unitary orbits, Matrix geometric mean.
\vskip 5pt\noindent
{\it 2010 mathematics subject classification.} 47A30, 15A60.
}

\section{Introduction } 

Let $\bM_n$ denote the space of complex $n\times n$ matrices and let $\bM_n^+$ stand for the positive (semi-definite) cone.
A linear map $\Phi:\bM_n\to\bM_m$ is called positive if $\Phi(\bM_n^+)\subset \bM_m^+$. A fundamental fact, the Russo-Dye Theorem \cite{RD}, asserts that every positive linear map attains its norm at the identity. This is discussed  in two nice  books, \cite[pp.\ 41--44]{Bhatia} and, in the operator algebras setting, \cite[Corollary 2.9]{Paulsen}. If $I$ denotes the identity and $\|\cdot\|_{\infty}$
 the operator norm, the theorem can be stated as:

\vskip 5pt 
\begin{theorem}\label{thRD}  If $\Phi: \bM_n\to \bM_m$ is  a positive linear map then,  for all contractions $A\in\bM_n$,
$$
\| \Phi(A)  \|_{\infty} \le \| \Phi(I) \|_{\infty}.
$$
\end{theorem}

\vskip 5pt
This short note presents several refinements of  Theorem \ref{thRD}. The contractive property of $A$ can be written as an operator inequality, $|A|\le I$ where $|A|=(A^*A)^{1/2}$ is the right modulus. This suggests that for a general matrix $Z$ the data of its two modulus $|Z|$ and $|Z^*|$ should be useful to obtain more general forms of the Russo-Dye Theorem. For instance, we have the following consequence of our main result.

\vskip 5pt
\begin{prop}\label{propintro} Let $Z\in\bM_n$, $J\in\bM_n^+$, and let $\Phi: \bM_n\to \bM_m$ be  a positive linear map. If $|Z| \le J$ and $|Z^*| \le J$,  then there exists  a unitary $V\in\bM_m$ such that
$$
 |\Phi(Z)|  \le \frac{\Phi(J) +V\Phi(J)V^*}{2}.
$$
\end{prop}

\vskip 5pt
If $Z=A$ is a contraction and $J=I$, Proposition \ref{propintro} claims that
\begin{equation}\label{eqintro}
 |\Phi(A)|  \le \frac{\Phi(I) +V\Phi(I)V^*}{2}
\end{equation}
and this clearly improves  Theorem \ref{thRD}. Let us illustrate that with some eigenvalue inequalities for the Schur product with $S\in\bM_n^+$. Let $s_{i}^{\downarrow}$ be the diagonal entries of $S$ rearranged in decreasing order, and similary let $\lambda^{\downarrow}_i(H)$ denote the eigenvalues of the Hermitian matrix $H$ arranged in  decreasing order. Theorem \ref{thRD} for the Schur multiplier $\Phi(T)=S\circ T$ says that
$$
\lambda^{\downarrow}_1(|S\circ A|) \le s_{1}^{\downarrow}
$$
for all contractions $A\in\bM_n$. From \eqref{eqintro} we infer several other inequalities, for example,
\begin{equation}\label{eqschur}
\lambda^{\downarrow}_3(|S\circ A|) \le s_{2}^{\downarrow}.
\end{equation}

 We shall prove a considerable improvement of Proposition \ref{propintro}. This yields several eigenvalue inequalities which extend and complete the Russo-Dye Theorem. Of course  the Russo-Dye Theorem holds for positive linear maps acting on unital $C^*$-algebras and,  at the end of the paper,  we state our main result in this setting.

\section{Positive linear maps and geometric mean}

\vskip 5pt
Let $Z\in\bM_n$, $J\in\bM_n^+$. We wish to state a version of Proposition \ref{propintro} under a more general assumption than the domination $|Z|, |Z^*|\le J$. Given two nonnegative functions $f(t)$ and $g(t)$, we use the notation
$$
Z\le_{f,g} J \quad \iff  \quad f(|Z|) \le J \  {\mathrm{and}}\  g(|Z^*|) \le J. 
$$
So, the assumption of Proposition \ref{propintro} can be written  $Z\le_{f,g} J$ with $f(t)=g(t)=t$.

We  also sharpen  Proposition \ref{propintro} by using the geometric mean instead of the arithmetic mean. We refer to  \cite{Ando} or \cite{Bhatia} for a background on the geometric mean $A\# B$ of two matrices $A,B\in\bM_n^+$. Since the arithmetic-geometric mean inequality
$$
A\#B \le \frac{A+B}{2}
$$
holds, the following theorem extends and improves Proposition \ref{propintro}.

\vskip 5pt
\begin{theorem}\label{th1} Let $Z\in\bM_n$, $J\in\bM_n^+$, and let $\Phi: \bM_n\to \bM_m$ be  a positive linear map. Suppose that $Z\le_{f,g} J$ for some   functions satisfying $f(t)g(t)=t^2$. Then, there exists  a unitary $V\in\bM_m$ such that
$$
 |\Phi(Z)|  \le \Phi(J)\#V\Phi(J)V^*.
$$
\end{theorem}

\vskip 5pt
\begin{proof} We first prove  the  case when $Z=N$  is normal, $f(t)=g(t)=t$, and  $J=|N|$. The theorem then reads as
\begin{equation}\label{eqt1}
 |\Phi(N)|  \le \Phi(|N|)\#V\Phi(|N|)V^*.
\end{equation}
Observe that
\begin{equation*}
\begin{pmatrix}
|N|& N^*\\
N& |N|
\end{pmatrix} \ge 0.
\end{equation*}
We may restrict $\Phi$ to the unital $C^*$-algebra spanned by $N$, and thus we may assume that $\Phi$ is completely positive  thanks to Stinespring's lemma  \cite{S}. Hence,
\begin{equation*}
\begin{pmatrix}
\Phi(|N|)& \Phi(N^*)\\
\Phi(N)& \Phi(|N|)
\end{pmatrix} \ge 0.
\end{equation*}
Next,  applying a unitary congruence 
\begin{equation*}
\begin{pmatrix}
I& 0\\
0& V
\end{pmatrix}
\begin{pmatrix}
\Phi(|N|)& \Phi(N^*)\\
\Phi(N)& \Phi(|N|)
\end{pmatrix}
\begin{pmatrix}
I& 0\\
0& V^*
\end{pmatrix}
\end{equation*}
where $V^*$ is the unitary  factor in the polar decomposition of $\Phi(N)$, we obtain
\begin{equation*}
\begin{pmatrix}
\Phi(|N|)& |\Phi(N)|\\
|\Phi(N)|& V\Phi(|N|)V^*
\end{pmatrix} \ge 0.
\end{equation*}
It then follows from the maximal property of $\#$ that \eqref{eqt1} holds.

Now, suppose that $Z=A$ is a contraction. Since $$|A|=\frac{W+W^*}{2}$$ for some unitary $W$, we infer from the polar decomposition that 
$$A=\frac{U_0+U_1}{2}$$some unitaries $U_0$ and $U_1$. Applying \eqref{eqt1} to the normal (unitary) matrix
\begin{equation*}
N=\begin{pmatrix}
U_0& 0\\
0& U_1
\end{pmatrix} 
\end{equation*}
and to the positive linear map from $\bM_{2}(\bM_n)$ to $\bM_n$,
\begin{equation*}
\begin{pmatrix}
P& Q\\
R& S
\end{pmatrix} \mapsto \Phi\left(\frac{P+S}{2}\right),
\end{equation*}
we obtain
\begin{equation}\label{eqt2}
   |\Phi(A)| \le \Phi(I) \#V\Phi(I)V^*.
\end{equation}
This establishes the theorem when $Z$ is a  contraction, $f(t)=g(t)=t$, and $J=I$.

We turn to the general case,  $f(|Z|)\le J$ and $g(|Z^*|)\le J$ with $f(t)g(t)=t^2$. Define a map $\Psi:\bM_n\to\bM_m$ by 
$$
\Psi(X) =\Phi(J^{1/2}XJ^{1/2}).
$$
Observe that  
\begin{equation}\label{eqt3}
\Psi(I)=\Phi(J), \quad {\mathrm{and}} \quad \Psi(Y)=\Phi(Z) 
\end{equation}
where 
$Y=J^{-1/2}ZJ^{-1/2},
$
and $J^{-1}$ is the generalized inverse. Thanks to the polar decomposition $Z=|Z^*|U=U|Z|$ and the condition $\sqrt{f(t)g(t)}=t$, 
$$Y=J^{-1/2}g^{1/2}(|Z^*|)Uf^{1/2}(|Z|)J^{-1/2}.
$$
From the assumption $Z\le_{f,g} J$ we infer that $J^{-1/2}g^{1/2}(|Z^*|)$ and $f^{1/2}(|Z|)J^{-1/2}$ are contractions, and so $Y$ is a contraction too. Applying \eqref{eqt2} to $Y$ and $\Psi$ yields
\begin{equation*}
   |\Psi(Y)| \le \Psi(I) \#V\Psi(I)V^*
\end{equation*}
and coming back to \eqref{eqt3} completes the proof.
\end{proof}

\vskip 5pt
In the course of the proof, we have noted two special cases.

\vskip 5pt
\begin{cor}\label{cornormal} Let $\Phi: \bM_n\to \bM_m$ be  a positive linear map  and let $N\in \bM_n$ be normal. Then, there exists  a unitary $V\in\bM_m$ such that
$$
 |\Phi(N)|  \le \Phi(|N|)\#V\Phi(|N|)V^*.
$$
\end{cor}

\vskip 5pt
\begin{cor}\label{corcontraction}  Let $A\in \bM_n$ be a contraction and  let $\Phi: \bM_n\to \bM_m$ be  a positive linear map. Then, for some unitary $V\in\bM_m$,
$$
   |\Phi(A)| \le \Phi(I) \#V\Phi(I)V^*.
$$
\end{cor}

\vskip 5pt
Another extension of the Russo-Dye Theorem is given in the next corollary.
\vskip 5pt
\begin{cor}\label{corrange}  Let $Z\in\bM_n$, $J\in\bM_n^+$, and let $\Phi: \bM_n\to \bM_m$ be  a positive linear map. If $J\ge E$, the range projection of $Z$, and $J\ge Z^*Z$, then, for some unitary $V\in\bM_m$,
$$
 |\Phi(Z)|  \le \Phi(J)\#V\Phi(J)V^*.
$$
\end{cor}

\vskip 5pt
\begin{proof} Apply Theorem \ref{th1} with $f(t)=t^2$ and $g(t)$ such that $g(0)=0$ and $g(s)=1$ for $s>0$.
\end{proof}

\vskip 5pt
The next statement extends Corollary \ref{cornormal}.

\vskip 5pt
\begin{cor}\label{cork}  Let $\Phi: \bM_n\to \bM_m$ be  a positive linear map, let $Z\in\bM_n$ be invertible, and  $\rho$  the spectral radius of $|Z^*| |Z|^{-1}$. Then, there exists  a unitary $V\in\bM_m$ such that
$$
 |\Phi(Z)|  \le \sqrt{\rho}\,\Phi(|Z|) \#V\Phi(|Z|)V^*.
$$
\end{cor}

\vskip 5pt
\begin{proof} Note that $|Z^*|\le \rho |Z|$ and apply Theorem \ref{th1}
for the matrices $Z$ and $J=\sqrt{\rho}|Z|$ with the functions $f(t)=\sqrt{\rho}t$ and $g(t)=(1/\sqrt{\rho})t$.
\end{proof}

\vskip 5pt
 In general the inequality of Corollary 2.5 is sharp as shown with the transpose map $\Phi(X)=X^T$ on $\bM_2$ and
$$
Z=\begin{pmatrix}
0&1 \\
k&0
\end{pmatrix}.
$$
Indeed, the spectral radius of $|Z^*||Z|^{-1}$ then equals $k$, and if $e_2$ denotes the second vector of the canonical basis of $\bC^2$, the inequality
$$
k=\langle e_2,|Z^T|e_2 \rangle \le c\langle e_2, |Z|^T\#V|Z|^TV^* e_2\rangle
$$
entails 
$$
k\le c\langle e_2, |Z|^Te_2\rangle \# \langle e_2, V|Z|^TV^* e_2\rangle  \le c\sqrt{k}
$$
and so $c\ge \sqrt{k}$.

For $2$-positive maps, Theorem \ref{th1} can be improved.

\vskip 5pt
\begin{theorem}\label{th2} Let $Z\in\bM_n$ and let $\Phi: \bM_n\to \bM_m$ be  a $2$-positive linear map. Then, for any pair of nonnegative functions $f(t)$ and $g(t)$    satisfying $f(t)g(t)=t^2$, there exists  a unitary $V\in\bM_m$ such that
$$
 |\Phi(Z)|  \le \Phi(f(|Z|))\#V\Phi(g(|Z^*|))V^*.
$$
\end{theorem}

\vskip 5pt
\begin{proof} From the polar decomposition $Z=U|Z|$ we have
$
Z= \sqrt{g(|Z^*|)} U \sqrt{f(|Z|)}$, hence
\begin{equation*}
\begin{pmatrix}
f(|Z|)& Z^*\\
Z& g(|Z^*|)
\end{pmatrix} \ge 0.
\end{equation*}
Since $\Phi$ is $2$-positive, we then have
\begin{equation*}
\begin{pmatrix}
\Phi(f(|Z|))& \Phi(Z^*)\\
\Phi(Z)& \Phi(g(|Z^*|))
\end{pmatrix} \ge 0.
\end{equation*}
Arguing as in the proof of Theorem \ref{th1} we then infer that
\begin{equation*}
 |\Phi(Z)|  \le \Phi(f(|Z|)) \#V\Phi(g(|Z^*|))V^*
\end{equation*}
for some unitary $V\in\bM_m$.
\end{proof}

\vskip 5pt
\begin{cor}\label{cor2pos}  Let $\Phi: \bM_n\to \bM_m$ be  a $2$-positive linear map, let $Z\in\bM_n$ and let $p\in(-\infty,\infty)$. Then, there exists  a unitary $V\in\bM_m$ such that
$$
 |\Phi(Z)|  \le \Phi(|Z|^{1+p}) \#V\Phi(|Z^*|^{1-p})V^*.
$$
\end{cor}

\vskip 5pt
Here, if $Z$ is not invertible, the nonpositive powers of $|Z|$ and $|Z^*|$ are understood as generalized inverses, in particular $|Z|^0$ is the support projection of $Z$ and $|Z^*|^0$ the range projection.
The case $p=0$ reads as
\begin{equation}\label{cex}
 |\Phi(Z)|  \le \Phi(|Z|) \#V\Phi(|Z^*|)V^*.
\end{equation}
However, in general \eqref{cex} does not hold if the 2-positivity assumption is dropped, as shown by the next example.

\vskip 5pt
\begin{example} Let $\Phi:\bM_2\to\bM_2$, $\Phi(X)=X+X^T$, and let 
$$Z=\begin{pmatrix}
0& 4\\
1& 0
\end{pmatrix} $$
Then, for any pair of unitary matrices $U,V\in\bM_2$,
$$
25=\det |\Phi(Z)| > \det \left\{U\Phi(|Z|)U^*\#  V\Phi(|Z^*|)V^*\right\}=16.
$$
\end{example}

\section{Some consequences}

\vskip 5pt
From these results follow some eigenvalue inequalities. Given $S,T\in\bM_n^+$, recall that
the (weak) log-majorisation $S \prec_{w\!\log}T$
means that a series of $n$ eigenvalue inequalities holds: For $k=1,\ldots,n$,
$$
\prod_{j=1}^k \lambda_j^{\downarrow}(S)
\le
\prod_{j=1}^k \lambda_j^{\downarrow}(T),
$$ 
where $\lambda_{j}^{\downarrow}(\cdot)$ still stand for the eigenvalues arranged in nonincreasing order.

Theorem \ref{th1} and Corollary \ref{corcontraction} are equivalent, though the most useful statement is Theorem \ref{th1}. For instance, if $Z\in\bM_n$ is a shift (a weighted permutation matrix), then  $|Z^*|$  and $|Z|$ are diagonal matrices, and comparing with a diagonal $J$ seems natural. A similar remark holds if $Z$ is paranormal. For a general matrix $Z$, the Kato supremum, \cite{Kato},
$$
|Z|\vee |Z^*| =\lim_{p\to\infty} \left(\frac{|Z|^p+|Z^*|^p}{2}\right)^{1/p}
$$
may be used as $J$ in Theorem \ref{th1} with $f(t)=g(t)=t$.  Some other natural choices for $J$ are $(|Z|^q+|Z^*|^q)^{1/q}$, $1\le q<\infty$, still with  $f(t)=g(t)=t$. 

We list  some consequences of Theorem \ref{th1}.

\vskip 5pt
\begin{cor}\label{corlog}   Let $Z\in\bM_n$, $J\in\bM_n^+$, and let $\Phi: \bM_n\to \bM_m$ be  a positive linear map. If $Z\le_{f,g} J$ for some   functions  satisfying $f(t)g(t)=t^2$, then, 
$$
   |\Phi(Z)| \prec_{w\!\log}  \Phi(J)
$$
and
$$
\lambda_{j+k+1}^{\downarrow } ( |\Phi(Z)|) \le \sqrt{\lambda_{j+1}^{\downarrow} ( \Phi(J)) \lambda_{k+1}^{\downarrow} ( \Phi(J))}
$$
 for all $j,k=0,1,\ldots,n$.
\end{cor}

\vskip 5pt
\begin{proof} By a well-known property of the geometric mean, Theorem \ref{th1} says that
$$
 |\Phi(Z)|  \le \sqrt{\Phi(J)}U \sqrt{\Phi(J)}V^*
$$
for some unitary matrices $U$ and $V$. From Horn's inequalities we then get the log-majorisation, while the second series of inequalities follows from a standard min-max principle.
\end{proof}

\vskip 5pt The Schur product of a matrix with the identity yields the diagonal part of the matrix. From the second statement of Corollary \ref{corlog} we obtain:

\vskip 5pt
\begin{cor}\label{corschur2}  Let $S\in\bM_n^+$ and let $A\in \bM_n$ be a contraction.  If $s_1^{\downarrow}\ge \ldots \ge s_n^{\downarrow}$ denote the diagonal entries of $S$ arranged in decreasing order, then,
for all  $j,k=0,1,\ldots$,
$$
 \lambda^{\downarrow}_{j+k+1}(  |S\circ A| ) \le \sqrt{s^{\downarrow}_{j+1}s^{\downarrow}_{k+1}}.
$$
\end{cor}

\vskip 5pt
\begin{remark}\label{corfisher}  For $j=k=1$, we recapture \eqref{eqschur}. If $S\in\bM_n^+$ is expansive,  Corollary \ref{corschur2} shows that
$$
 \lambda^{\downarrow}_{2j+1}(  S\circ S^{-1} ) \le s^{\downarrow}_{j+1}. 
$$
  If $S\in\bM_n^+$ is a contraction, then
$$
 \lambda^{\downarrow}_{2j+1}(  S\circ S ) \le s^{\downarrow}_{j+1}. 
$$
\end{remark}

\vskip 5pt
\begin{cor}\label{corprod} Let $Z\in\bM_n$, $J\in\bM_n^+$, and let $\Phi: \bM_n\to \bM_m$ be  a positive linear map. If $Z\le_{f,g} J$ for some   functions  satisfying $f(t)g(t)=t^2$, then
$$
 \left\{\prod_{j=1}^k \lambda_j^{\uparrow} ( |\Phi(Z)|)\right\}^2 \le \prod_{j=1}^k \lambda_j^{\uparrow}(\Phi(J) \lambda_j^{\downarrow}(\Phi(J))
$$
for all $k=1,\ldots, n$.
\end{cor}

\vskip 5pt
\begin{proof} Let ${\mathcal{S}}$ be any subspace of $\bC^m$, and denote by $A_{\mathcal{S}}$ the compression of $A\in\bM_m$ onto ${\mathcal{S}}$. From Theorem \ref{th1} and Ando's  inequality \cite[Theorem 3]{Ando} we infer
$$
 |\Phi(Z)|_{\mathcal{S}}  \le \left(\Phi(J)\right)_{\mathcal{S}}\#\left(V\Phi(J)V^*\right)_{\mathcal{S}}.
$$
So, 
$$
\det |\Phi(Z)|_{\mathcal{S}}  \le \left\{\det \left(\Phi(J)\right)_{\mathcal{S}}\det \left(\Phi(J)\right)_{V^*(\mathcal{S})} \right\}^{1/2}.
$$
Letting ${\mathcal{S}}$ be the spectral subspace corresponding to the $k$ smallest eigenvalues of $\Phi(J)$ and using standard min-max principles completes the proof.
\end{proof}

Every  $n$-by-$n$ complex matrix $Z$ has a decomposition in real and imaginary parts, its so-called Cartesian decomposition, $Z=X+iY$ where $X$ and $Y$ are self-adjoints. In general, the operator inequality 
\begin{equation}\label{eq0}|Z|\le |X|+|Y|
\end{equation}
 does not hold. However, a well-known substitute is the (weak) log-majorisation
\begin{equation}\label{eq1}
   |Z| \prec_{w\!\log}  |X|+|Y|,
\end{equation}
see \cite[Corollary 2.11]{BL1} for a stronger statement. We may use Corollary \ref{cornormal} to obtain further extensions of the triangle inequality \eqref{eq1}.

\vskip 5pt
\begin{cor}\label{corcart} Let $\Phi: \bM_n\to \bM_m$ be  a positive linear map  and let $Z\in\bM_n$ with Cartesian decomposition $Z=X+iY$. Then, there exists  a unitary $V\in\bM_m$ such that
$$
 |\Phi(Z)|  \le \Phi(|X|+|Y|)\#V\Phi(|X|+|Y|)V^*.
$$
\end{cor}

\vskip 5pt
\begin{proof}
Consider the normal matrix in $\bM_{2n}$
\begin{equation*}N=
\begin{pmatrix}
X& 0\\
0& iY
\end{pmatrix}
\end{equation*}
and apply Corollary \ref{cornormal} to this matrix and the positive linear map $\Psi:\bM_{2n}\to \bM_m$ defined as $\Psi=\Phi \circ \Lambda$ where $\Lambda$ is the partial trace on $\bM_{2n}=\bM_2(\bM_n)$, 
\begin{equation*}\Lambda\left(\begin{pmatrix}
A& B\\
C& D
\end{pmatrix}\right) = A+D.
\end{equation*}
Since $\Psi(N)=\Phi(Z)$ and $\Psi(|N|)=\Phi(|X|+|Y|)$, the corollary is established.
\end{proof}

\vskip 5pt
\begin{cor}\label{corcartlog}  If $\Phi: \bM_n\to \bM_m$ is  a positive linear map  and  $Z\in\bM_n$ has Cartesian decomposition $Z=X+iY$, then, 
$$
   |\Phi(Z)| \prec_{w\!\log}  \Phi(|X|+|Y|).
$$
\end{cor}

\vskip 5pt
\begin{cor}\label{cornorm}  If   $Z\in\bM_n$ has Cartesian decomposition $Z=X+iY$, then, 
$$
   \left\| (|X|+|Y|)^{-1/2} Z  (|X|+|Y|)^{-1/2}\right\|_{\infty} \le 1.
$$
\end{cor}

\vskip 5pt
\begin{proof} This is a direct consequence of Corollary \ref{corcartlog} applied to the positive linear map $T\mapsto (|X|+|Y|)^{-1/2} T  (|X|+|Y|)^{-1/2}$.
\end{proof}

\vskip 5pt
\begin{remark} Corollary \ref{cornorm} is surprising. We could think that this follows from the polar decomposition $Z=|Z^*|^{1/2}U|Z|^{1/2}$ and the plausible fact that $|Z|^{1/2}(|X|+|Y|)^{-1/2}$ is contractive. But in general this matrix is not contractive, as it would be equivalent to the untrue inequality \eqref{eq0}.
\end{remark}

\vskip 5pt
\begin{cor}\label{corradius}  If   $Z\in\bM_n$ has Cartesian decomposition $Z=X+iY$, then, 
$$
   \rho( Z  (|X|+|Y|)^{-1}) \le 1.
$$
\end{cor}

\vskip 5pt
\begin{proof} This is a straightforward consequence of Corollary \ref{cornorm}.
\end{proof}

\vskip 5pt
\begin{remark} Corollary \ref{corradius} seems surprising too: It may happen that 
$$
   \left\| Z  (|X|+|Y|)^{-1})\right\|_{\infty} > 1.
$$ 
Otherwise we would have $|Z|^2\le  (|X|+|Y|)^2 $, entailing the untrue inequality \eqref{eq0}.
\end{remark}

\section{Comments}

This article is a continuation of  \cite{BL2}. Indeed, 
the idea of revisiting the Russo-Dye Theorem comes from    that paper, where Corollary \ref{cornormal} is proved with several related sharp inequalities. 
 We can state  a version of Theorem \ref{th1} for operator algebras. Let $\mathcal{A}$ be a unital $C^*$-algebra, and let $\mathcal{B}$ be a von Neumann algebra.

\vskip 5pt
\begin{theorem}\label{th2} Let $Z\in{\mathcal{A}}$, $J\in{\mathcal{A}}^+$, and let $\Phi: {\mathcal{A}}\to {\mathcal{B}}$ be  a positive linear map. Suppose that $Z\le_{f,g} J$ for some  continuous functions satisfying $f(t)g(t)=t^2$. Then, there exists  a partial isometry $V\in{\mathcal{B}}$ such that
$$
 |\Phi(Z)|  \le \Phi(J)\#V\Phi(J)V^*.
$$
\end{theorem}

\vskip 5pt
The proof follows exactly the same steps as in the matrix case. Note that $V^*$ is the partial isometry occuring in the polar decomposition $\Phi(Z)=V^*|\Phi(Z)|$, in general it is not a unitary operator, and it belongs to the von Neumann algebra generated by $\Phi(Z)$. If ${\mathcal{A}}$ is also a von Neumann algebra, then $f(t)$ and $g(t)$ can be bounded Borel functions.

The assumption of Corollary \ref{cork} can be replaced by the condition that $|Z^*|\le \rho|Z|$ for some $\rho>0$, which makes sense also for noninvertible operators. With $\rho =1$, Corollary \ref{cork} then takes  the following  form for infinite dimensional operators.

\vskip 5pt
\begin{cor}\label{corsemi} Let $Z\in{\mathcal{A}}$ be semi-hyponormal and let $\Phi: {\mathcal{A}}\to {\mathcal{B}}$ be  a positive linear map. Then, there exists  a partial isometry $V\in{\mathcal{B}}$ such that
$$
 |\Phi(Z)|  \le \Phi(|Z|)\#V\Phi(|Z|)V^*.
$$
\end{cor}

\vskip 5pt
\noindent
Laboratoire de math\'ematiques, 

\noindent
Universit\'e de Bourgogne Franche-Comt\'e, 

\noindent
25 000 Besan\c{c}on, France.

\noindent
Email: jcbourin@univ-fcomte.fr

  \vskip 10pt
\noindent Department of mathematics, KNU-Center for Nonlinear
Dynamics,

\noindent
Kyungpook National University,

\noindent
 Daegu 702-701, Korea.

\noindent Email: eylee89@knu.ac.kr

\end{document}